\documentclass{amsart}
\usepackage{amsfonts}

\usepackage{graphicx}
\usepackage{amscd}

\numberwithin{equation}{section}
\newtheorem{theorem}{Theorem}[section]
\theoremstyle{plain}

\newtheorem{corollary}[theorem]{Corollary}

\newtheorem{lemma}[theorem]{Lemma}

\newtheorem{proposition}[theorem]{Proposition}

\numberwithin{equation}{section}

\begin{document}
\title[Twisted Fock spaces ]{Fock spaces corresponding to positive definite linear transformations}
\author{R. Fabec, G. \'{O}lafsson, and A. N. Sengupta}
\address{Department of mathematics\\
Louisiana State University\\
Baton Rouge\\
LA 70803} \email{fabec@math.lsu.edu}
\urladdr{http://www.math.lsu.edu/\symbol{126}fabec}
\email{olafsson@math.lsu.edu}
\urladdr{http://www.math.lsu.edu/\symbol{126}olafsson}
\email{sengupta@math.lsu.edu}
\urladdr{http://www.math.lsu.edu/\symbol{126}sengupta}
\thanks{The research of G. \'Olafsson was supported by DMS-0070607 and DMS-0139473}
\thanks{The research of A. Sengupta was supported by DMS-0201683}

\begin{abstract}
Suppose $A$ is a positive real linear transformation on a finite
dimensional complex inner product space $V$. The reproducing
kernel for the Fock space of square integrable holomorphic
functions on $V$ relative to the Gaussian measure $d\mu_A(z)=\frac
{\sqrt {\det A}} {\pi^n}e^{-\mbox{\rm Re}\langle Az,z\rangle}\,dz$
is described in terms of the holomorphic--antiholomorphic
decomposition of the linear operator $A$. Moreover, if $A$
commutes with a conjugation on $V$, then a restriction mapping to
the real vectors in $V$ is polarized to obtain a Segal--Bargmann
transform, which we also study in the Gaussian-measure setting.
\end{abstract}

\maketitle

\section*{Introduction}

The classical Segal-Bargmann transform is an integral transform
which defines a  unitary isomorphism of $L^{2}(\mathbb{R}^{n})$
onto the Hilbert space $\mathbf{F}(\mathbb{C}^{n})$ of entire
functions on $\mathbb{C}^{n}$ which are square integrable with
respect to the Gaussian measure $\mu =\pi
^{-n}e^{-|z|^{2}}dxdy$, where $dxdy$ stands for the Lebesgue measure on $%
\mathbb{R}^{2n}\simeq \mathbb{C}^{n}$, see \cite{B61,F28,F32,Fo89,S27,S63}.
There have been several generalizations of this transform, based on the heat
equation or the representation theory of Lie groups \cite{H94,GOO96,AS99}.
In particular, it was shown in \cite{GOO96} that the Segal-Bargmann transform
is a special case of \textit{the restriction principle}, i.e., construction
of unitary isomorphisms based on the polarization of a restriction map. This
principle was first introduced in \cite{GOO96}, see also \cite{GO00}, where
several examples were explained from that point of view.\ In short the
restriction principle can be explained in the following way. Let $M_{\mathbb{%
C}}$ be a complex manifold and let $M\subset M_{\mathbb{C}}$ be a totally
real submanifold. Let $\mathbf{F}=\mathbf{F}(M_{\mathbb{C}})$ be a Hilbert
space of holomorphic functions on $M_{\mathbb{C}}$ such that the evaluation
maps $\mathbf{F}\ni F\mapsto F(z)\in \mathbb{C}$ are continuous for all $%
z\in M_{\mathbb{C}}$, i.e., $\mathbf{F}$ is a \textit{reproducing Hilbert
space}. There exists a function $K:M_{\mathbb{C}}\times M_{\mathbb{C}%
}\rightarrow \mathbb{C}$ holomorphic in the first variable, anti-holomorphic
in the second variable, and such that the following hold:

\begin{enumerate}
\item  $K(z,w)=\overline{K(w,z)}$ for all $z,w\in M_{\mathbb{C}}$;

\item  If $K_{w}(z):=K(z,w)$ then $K_{w}\in \mathbf{F}$ and
\begin{equation*}
F(w)=(F,K_{w}),\quad \forall F\in \mathbf{F},z\in M_{\mathbb{C}}\,.
\end{equation*}
\end{enumerate}

The function $K$ is the \textit{reproducing kernel} for the Hilbert space.
Let $D:M\rightarrow \mathbb{C}^{\ast }$ be measurable. Then the restriction
map $RF:=DF|_M$ is injective. Assume that there is a measure $\mu $ on $M$ such
that $RF\in L^{2}(M,\mu )$ for all $F$ in a dense subset of
$\mathbb F$.
Provided $R$ is closeable, polarizing $R^{\ast }$ we can write
\begin{equation*}
R^{\ast }=U|R^{\ast }|
\end{equation*}
where $U:L^{2}(M,\mu )\rightarrow \mathbf{F}$ is a  unitary
isomorphism. Using that $\mathbf{F}$ is a reproducing Hilbert
space we get that
\begin{equation*}
Uf(z)=(Uf,K_{z})=(f,U^{\ast }K_{z})=\int_{M}f(m)(U^{\ast }K_{z})(m)\,d\mu
(m)\,.
\end{equation*}
Thus $Uf$ is always an integral operator. We notice also that the formula
for $U$ shows that the important object in this analysis is the reproducing
kernel $K(z,w)$. The reproducing kernel for the classical Fock space is
given by $K(z,w)=e^{z\bar{w}}$. By taking $D(x):=(2\pi )^{-n/4}e^{-|x|^{2}}$%
, which is closely related to the heat kernel, we arrive at the classical
Segal--Bargmann transform
\begin{equation*}
Ug(x)=(2/\pi )^{n/4}e^{\left( x,x\right) /2}\int g(y)e^{-\left(
x-y,x-y\right) }\,dy\,.
\end{equation*}
The same principle can be used to construct the Hall--transform
for compact Lie groups, \cite{H94}. In \cite{DH99}, Driver and
Hall, motivated by application to quantum Yang-Mills theory,
introduced a Fock space and Segal--Bargmann transform depending on
two parameters $r,s>0$, giving
different weights to the $x$ and $y$ directions, where $z=x+iy\in \mathbb{C}%
^{n}$ (this was also studied  in \cite{AS99}). Thus $\mathbf{F}$
is now the space of holomorphic functions $F(z)$ on
$\mathbb{C}^{n}$ which are square-integrable with respect to the
Gaussian
measure $dM_{r,s}(z)=\frac{1}{(\pi r)^{n/2}(\pi s)^{n/2}}\,e^{-\frac{x^{2}}{r%
}-\frac{y^{2}}{s}}\,$. In \cite{AS99} the reproducing kernel and the
Segal--Bargmann transform for this space is worked out. This construction
has a natural generalization by viewing $r^{-1}$ and $s^{-1}$ as the
diagonal elements in a positive definite matrix $A=d(r^{-1}I_{n},s^{-1}I_{n})
$. The measure is then simply
\begin{equation}
dM_{r,s}(z)=\frac{\sqrt{\det (A)}}{\pi ^{n}}e^{-(Az,z)}\,dxdy  \label{eq1}
\end{equation}
and this has meaning for any positive definite matrix $A$.

In this paper we show that (\ref{eq1}) gives rise to a Fock space $\mathbf{F}%
_{A}$ for arbitrary positive matrices $A$. We find an expression
for the reproducing kernel $K_{A}(z,w)$. We use the restriction
principle to construct a natural generalization of the
Segal-Bargmann transform for this space, with a   certain natural
restriction on $A$. We study this also in the Gaussian setting,
and indicate a generalization to infinite dimensions.

\section{The Fock space and the restriction principle}

In this section we recall some standard facts about the classical Fock space
of holomorphic function on $\mathbb{C}^{n}$. We refer to \cite{Fo89} for
details and further information. Let $\mu $ be the measure $d\mu =\pi
^{-n}e^{-\left| |z\right| |^{2}}dxdy$ and let $\mathbf{F}$ be the classical
Fock-space of holomorphic functions $F:\mathbb{C}^{n}\rightarrow \mathbb{C}$
such that
\begin{equation*}
\left| \left| F\right| \right| ^{2}:=\int \left| F(z)\right| ^{2}\,d\mu(z) <\infty \,.
\end{equation*}
(Note that the term ``Fock space'' is also used for the completed
symmetric tensor algebra over a Hilbert space, but that is not our
usage here.) The space $\mathbf{F}$ is a reproducing Hilbert space
with inner product
\begin{equation*}
\left( F,G\right) =\int F(z)\overline{G(z)}\,\,d\mu \,
\end{equation*}
and reproducing kernel $K(z,w)=e^{\langle z,w\rangle}$, where
$\langle z,w\rangle=z{\overline w}=z_1{\overline w_1}+\cdots
z_n{\overline w_n}$. Thus
\begin{equation*}
F(w)=\int F(z)\overline{K(z,w)}\,d\mu =(F,K_{w})
\end{equation*}
where $K_{w}(z)=K(z,w)$. The function $K(z,w)$ is holomorphic in the first
variable, anti-holomorphic in the second variable, and $K(z,w)=\overline{%
K(w,z)}$. Notice that $K(z,z)=(K_{z},K_{z})$. Hence $||K_{z}||=e^{|z|^{2}/2}$%
. Finally the linear space of finite linear combinations $\sum
c_{j}K_{z_{j}} $, $z_{j}\in \mathbb{C}^{n}$, $c_{j}\in \mathbb{C}$, is dense
in $\mathbf{F}$. An orthonormal system in $\mathbf{F}$ is given by the
monomials $e_{\alpha }(z)=z_{1}^{\alpha _{1}}\cdots z_{n}^{\alpha _{n}}/%
\sqrt{\alpha _{1}!\cdots \alpha _{n}!}$, $\alpha \in \mathbb{N}_{o}^{n}$.

View $\mathbb{R}^{n}\subset \mathbb{C}^{n}$ as a totally real submanifold of
$\mathbb{C}^{n}$. We will now recall the construction of the classical
Segal-Bargmann transform using the \textit{restriction principle}, see \cite
{GO00,GOO96}.
For constructing a restriction map as explained in the introduction
we need to choose the function $D(x)$.  One motivation for the choice
of $D$ is the heat kernel, but another one, more closely related to
representation theory, is that the restriction map should commute
with the action of $\mathbb{R}^n$ on the Fock space and
$L^2(\mathbb{R}^n)$.  Indeed, take
\begin{equation*}
T(x)F(z)=m(x,z)F(z-x)
\end{equation*}
for $F$ in $\mathbf{F}$ where $m(x,z)$ has properties sufficient to make
$x\mapsto T(x)$ a unitary representation of $\mathbb{R}^n$ on $\mathbf{F}$.
Namely, $m$ is a
multiplier, i.e.,  $m(x,z)m(y,z-x)=m(x+y,z)$; $z\mapsto m(x,z)$ is
holomorphic in $z$ for each $x$; and
$|m(x,z)|=\left(\frac {d\mu (z-x)}{d\mu (z)}%
\right)^{\frac 12}=e^{(\text{\textrm{Re}}z\cdot x-||x||^2/2)}$.  Note
$m(x,z):=e^{z\cdot x-||x||^2/2}$ has these properties.  Set
$D(x)=(2\pi)^{-n/4}m(0,x)=(2\pi )^{-n/4}e^{-||x||^2/2}$ and define
$R:\mathbf{F}\rightarrow C^{\infty}(%
\mathbb{R}^n)$ by
\begin{equation*}
RF(x):=D(x)F(x)=(2\pi )^{-n/4}e^{-||x||^{2}/2}F(x).
\end{equation*}

Then
\begin{align*}
RT(y)F(x)& =(2\pi )^{-n/4}e^{-||x||^{2}/2}T(y)F(x) \\
& =(2\pi )^{-n/4}e^{-||x||^{2}/2}e^{x\cdot y-||y||^{2}/2}F(x-y) \\
& =(2\pi )^{-n/4}e^{-||x-y||^{2}/2}F(x-y) \\
& =RF(x-y).
\end{align*}
As $\mathbb{R}^{n}$ is a totally real submanifold of $\mathbb{C}^{n}$, it
follows that $R$ is injective. Furthermore the holomorphic polynomials $%
p(z)=\sum a_{\alpha }z^{\alpha }$ are dense in $\mathbf{F}$ and obviously $%
Rp\in L^{2}(\mathbb{R}^{n})$. Hence all the Hermite functions $h_{\alpha
}(x)=\left( -1\right) ^{\left| \alpha \right| }\left( D^{\alpha }e^{-||
x|| ^{2}}\right) e^{||x||^{2}/2}$ are in the image of $R$; so
$\mathrm{Im}(R)$ is dense in $\mathbf{L}^{2}(\mathbb{R}^{n})$ and $R$ is a
densely defined operator from $\mathbf{F}$ into $L^{2}(\mathbb{R}^{n})$. It
follows easily from the fact that the maps $F\mapsto F(z)$ are continuous,
that $R$ is a closed operator. \ Hence $R$ has an adjoint $R^{\ast }:L^{2}(%
\mathbb{R}^{n})\rightarrow \mathbf{F}$. For $z,w\in \mathbb{C}^{n}$, let $%
(z,w)=\sum z_{j}w_{j}$. Then:
\begin{eqnarray*}
R^{\ast }g(z) &=&(R^{\ast }g,K_{z}) \\
&=&(g,RK_{z}) \\
&=&(2\pi )^{-n/4}\int g(y)e^{-|\left| y\right| |^{2}/2}e^{z\cdot y}\,dy \\
&=&(2\pi )^{-n/4}e^{\left( z,z\right) /2}\int g(y)e^{-\left( z-y,z-y\right)
/2}\,dy \\
&=&(2\pi )^{n/4}e^{(z,z)/2}g\ast p(z)\,
\end{eqnarray*}
where $p(z)=(2\pi )^{-n/2}e^{-(z,z)/2}$ is holomorphic. Hence
\begin{equation}
RR^{\ast }g(x)=g\ast p(x)\,.  \label{eq3}
\end{equation}
As $p\in L^{1}(\mathbb{R}^{n})$, it follows that $\left| \left| RR^{\ast
}\right| \right| \leq \left| \left| p\right| \right|_{1}$; so $RR^{\ast }$
is continuous.
\begin{equation*}
(R^{\ast }g,R^{\ast }g)=(RR^{\ast }g,g)\leq \left| \left| RR^{\ast }\right|
\right| \,\left| \left| g\right| \right| _{2}\,.
\end{equation*}
Thus

\begin{lemma}
The maps $R$ and $R^{\ast }$ are continuous.
\end{lemma}

Let $p_{t}(x)=(2\pi t)^{-n/2}e^{-\left( x,x\right) /2t}$ be the heat kernel
on \ $\mathbb{R}^{n}$. Then $\left( p_{t}\right) _{t>0}$ is a convolution
semigroup and $p=p_{1}$. Hence $\sqrt{RR^{\ast }}=p_{1/2}\ast $ or
\begin{equation*}
RUg(x)=\left| R^{\ast }\right| g(x)=p_{1/2}\ast g(x)=\pi ^{-n/2}\int
g(y)e^{-\left( x-y,x-y\right) }\,dy\,.\,
\end{equation*}
It follows that
\begin{equation*}
Ug(x)=(2/\pi )^{n/4}e^{\left( x,x\right) /2}\int g(y)e^{-\left(
x-y,x-y\right) }\,dy
\end{equation*}
for $x\in \mathbb{R}^{n}$. But the function on the right hand side is
holomorphic in $x$. Analytic continuation gives the following theorem.

\begin{theorem}
The map $U:L^{2}(\mathbb{R}^{n})\rightarrow \mathbf{F}$ given by
\begin{equation*}
Ug(z)=(2/\pi )^{n/4}\int g(y)\exp (-\left( y,y\right) +2(z,y)-\left(
z,z\right) /2)\,dy
\end{equation*}
is a  unitary isomorphism. $U$ is called the Segal--Bargmann
transform.
\end{theorem}

\section{Twisted Fock spaces}

Let $V\simeq \mathbb{C}^{n}$ be a finite dimensional complex vector space of
complex dimension $n$ and let $\langle \cdot ,\cdot \rangle$ be a complex inner product.
As before we will sometimes write $\langle z,w\rangle=z\cdot w$.  We will also consider $%
V$ as a real vector space with real inner product defined by $(z,w)=\text{%
\textrm{Re}}\langle z,w\rangle$. Notice that $(z,z)=\langle z,z\rangle$
for all $z\in \mathbb{C}^{n}$.
Let $J$ be the real linear transformation of $V$ given by $Jz=iz$. Note that
$J^{\ast }=-J=J^{-1}$ and thus $J$ is a skew symmetric real linear
transformation. Fix a real linear transformation $A$. Then $A=H+K$ where
\begin{equation*}
H:=\frac{A+J^{-1}AJ}{2}\quad \text{and}\quad K:=\frac{A-J^{-1}AJ}{2}.
\end{equation*}
Note that $HJ=\frac{1}{2}(AJ-J^{-1}A)=\frac{1}{2}J(J^{-1}AJ+A)=JH$ and $KJ=%
\frac{1}{2}(AJ+J^{-1}A)=\frac{1}{2}J(J^{-1}AJ-A)=-JK$. Furthermore
$H$ is complex linear and $K$ is conjugate linear. We assume that
$A$ is symmetric and positive definite.

\begin{lemma}
The complex linear transformation $H$ is self adjoint, positive with
respect to the inner product $\langle \cdot ,\cdot \rangle$, and invertible.
\end{lemma}

\begin{proof}
Since $A$ is positive and invertible as a real linear transformation, we
have $(Az,z)>0$ for all $z\neq 0$. But $J$ is real linear and skew
symmetric. Hence $(JAJ^{-1}z,z)>0$ for all $z\neq 0$. In particular $H=\frac{%
1}{2}(A+JAJ^{-1})$ is complex linear, symmetric with respect to the real
inner product $(\cdot ,\cdot )$, and positive. We know $(Hv,w)=(v,Hw)$. Thus
$\text{\textrm{Re}}\langle Hv,w\rangle =\text{\textrm{Re}}\langle
v,Hw\rangle $. From this we obtain
\begin{equation*}
\text{\textrm{Re}}\langle Hiv,w\rangle =\text{\textrm{Re}}\langle
iv,Hw\rangle .
\end{equation*}
This implies $\text{\textrm{Im}}\langle Hv,w\rangle =\text{\textrm{Im}}%
\langle v,Hw\rangle $. Putting these together gives $\langle Hv,w\rangle
=\langle v,Hw\rangle $. Hence $H$ is complex self adjoint and $\langle
Hz,z\rangle >0$ for $z\neq 0$.
\end{proof}

\begin{lemma}
Let $w\in V$. Then $\langle Aw,w\rangle=(Aw,w)+i\mathrm{Im}\langle Kw,w\rangle$ and $%
(Aw,w)=(Hw,w)+(Kw,w)$.
\end{lemma}

\begin{proof}
Let $w\in V$. Then
\begin{align*}
\langle Aw,w\rangle & =\langle Hw,w\rangle +\langle Kw,w\rangle  \\
& =(Hw,w)+(Kw,w)+i\text{\textrm{Im}}\langle Kw,w\rangle  \\
& =(Aw,w)+i\text{\textrm{Im}}\langle Kw,w\rangle .
\end{align*}
This implies the first statement. Taking the real part in the
second line gives the second claim, which also follows directly
from bilinearity of $(\cdot,\cdot)$.
\end{proof}

Denote by $\det_{V}$ the determinant of a $\mathbb{R}$-linear map on $%
\mathbb{C}^{n}\simeq \mathbb{R}^{2n}$. Let $d\mu _{A}(z)=\pi ^{-n}\sqrt{%
\det_{V}A}e^{-(Az,z)}dxdy$ and let $\mathbf{F}_{A}$ be the space of
holomorphic functions $F:\mathbb{C}^{n}\rightarrow \mathbb{C}$ such that
\begin{equation*}
||F||_{A}^{2}:=\int |F(z)|^{2}\,d\mu _{A}<\infty \,.
\end{equation*}
Our normalization of $d\mu $ is chosen so that $||1||_{A}=1$. Just
as in the classical case one can show that $\mathbf{F}_{A}$ is a
reproducing Hilbert space, but this will also follow from the
following Lemma. We notice that all the holomorphic polynomials
$p(z)$ are in $\mathbf{F}$. To simplify the notation, we let
$T_1=H^{-1/2}$. Then $T_1$ is symmetric, positive definite and
complex linear. Let $c_{A}=\sqrt{\det_{V}(A^{1/2}T_1)}=\left(
\det_{V}(A)/\det_{V}(H)\right) ^{1/4}$.

\begin{lemma}
\label{isofock}Let $F:V\rightarrow \mathbb{C}$ be holomorphic.
Then $F\in \mathbf{F}_{A}$ if and only if $F\circ T_1\in
\mathbf{F}$ and the map $\Psi:\mathbf{F}\rightarrow
\mathbf{F}_{A}$ given by
\begin{equation*}
\Psi (F)(w):=c_{A}\exp \left( -\overline{\langle
KT_1w,T_1w\rangle}/2\right) F(T_1w)
\end{equation*}
is a  unitary isomorphism. In particular
\begin{equation*}
\Psi ^{\ast }F(w)=\Psi ^{-1}F(w)=c_{A}^{-1}\exp \left( \overline{\langle Kw,w\rangle}%
/2\right) F(\sqrt{H}w)\,.
\end{equation*}
\end{lemma}

\begin{proof}
Let $F:V\rightarrow \mathbb{C}$. Then $F$ is holomorphic if and
only if $F\circ T_1$ is holomorphic as $T_1$ is complex linear and
invertible. Moreover, we also have:
\begin{align*}
||\Psi F||^{2}& =\pi ^{-n}\int |\Psi F(w)|^{2}\,e^{-\langle w,w\rangle }\,dw
\\
& =\pi ^{-n}\sqrt{\det {}_{V}A}\int |F(w)|^{2}e^{-(Kw,w)}e^{-\langle\sqrt{H}w,%
\sqrt{H}w\rangle}\,dw \\
& =\pi ^{-n}\sqrt{\det {}_{V}A}\int |F(w)|^{2}e^{-(Kw,w)}e^{-\langle Hw,w\rangle}\,dw \\
& =\pi ^{-n}\sqrt{\det {}_{V}A}\int |F(w)|^{2}e^{-((H+K)w,w)}\,dw \\
& =\pi ^{-n}\sqrt{\det {}_{V}A}\int |F(w)|^{2}e^{-(Aw,w)}\,dw \\
& =||F||_{A}^{2}\,
\end{align*}
and thus, by polarization, $\Psi$ is unitary.
\end{proof}

\begin{theorem}
The space $\mathbf{F}_{A}$ is a reproducing Hilbert space with reproducing
kernel
\begin{equation*}
K_{A}(z,w)=c_{A}^{-2}e^{\frac{1}{2}\overline{\langle Kz,z\rangle }%
}e^{\langle Hz,w\rangle }e^{\frac{1}{2}\langle Kw,w\rangle }\,.
\end{equation*}
\end{theorem}

\begin{proof}
By Lemma \ref{isofock} we get
\begin{eqnarray*}
c_{A}\exp (-\overline{\langle KT_1w,T_1w\rangle}/2)F(T_1w) &=&\Psi (F)(w) \\
&=&(\Psi (F),K_{w})_{\mathbf{F}_{A}} \\
&=&(F,\Psi ^{\ast }(K_{w}))_{\mathbf{F}}\,.
\end{eqnarray*}
Hence
\begin{equation*}
K_{A}(z,w)=c_{A}^{-1}\exp (\overline{\langle Kw,w\rangle}/2)\Psi ^{\ast }(K_{\sqrt{H}%
w})=c_{A}^{-2}e^{\frac{1}{2}\overline{\langle Kz,z\rangle }}e^{\langle
Hz,w\rangle }e^{\frac{1}{2}\langle Kw,w\rangle }.
\end{equation*}
\end{proof}

\section{The Restriction Map}

We assume as before that $A>0$. We notice that Lemma \ref{isofock}
gives a unitary isomorphism $\Psi ^{\ast }U:L^{2}(\mathbb{R}^{n})
\rightarrow \mathbf{F}_{A}$, where $U$ is the classical
Segal-Bargmann transform. But this is not the natural transform that
we are looking for. As $H$ is positive definite there is an
orthonormal basis $e_{1},\ldots ,e_{n}$ of $V$ and
positive numbers $\lambda _{j}>0$ such that $He_{j}=\lambda _{j}e_{j}
$. Let $V_{\mathbb{R}}:=\sum \mathbb{R}e_{k}$. Set $\sigma (\sum a_
{i}e_{i})=\sum
\bar{a}_{i}e_{i}$. Then $\sigma$ is a conjugation with $V_{\mathbb{R}%
}=\{z:\sigma z=z\}$. We say that a vector is \textit{real} if it belongs to $%
V_{\mathbb{R}}$. As $He_{j}=\lambda _{j}e_{j}$ with $\lambda _{j}\in \mathbb{%
R}$ it follows that $HV_{\mathbb{R}}\subseteq V_{\mathbb{\ R}}$. We denote
by $\det $ the determinant of a $\mathbb{R}$-linear map of $V_{\mathbb{R}}$.

\begin{lemma}\label{l5}
$\langle Kz,w\rangle =\langle Kw,z\rangle $.
\end{lemma}

\begin{proof}
Note that $\sigma K$ is complex linear. Since $J^{\ast }=-J$,
$K=\frac{1}{2}(A-JAJ^{-1})$ is real symmetric. Thus $%
(Kw,z)=(w,Kz)=(Kz,w)$. Also note $(iKz,w)=(JKz,w)=-(KJz,w)=-(Jz,Kw)=-(iz,Kw)$%
. Hence $\text{\textrm{Re}}\langle iKz,w\rangle =-\text{\textrm{Re}}\langle
iz,Kw\rangle $. So $-\text{\textrm{Im}}\langle Kz,w\rangle =\text{\textrm{Im}%
}\langle z,Kw\rangle $. This gives $\text{\textrm{Im}}\langle Kw,z\rangle =%
\text{\textrm{Im}}\langle Kz,w\rangle $. Hence $\langle Kz,w\rangle =\langle
Kw,z\rangle $.
\end{proof}

\begin{lemma}
$(\sigma K)^{\ast }=K\sigma $.
\end{lemma}

\begin{proof}
We have $\langle \sigma z,\sigma w\rangle =\langle w,z\rangle $.
Hence\newline $\langle\sigma Kz,w\rangle=\langle\sigma w,\sigma
^{2}Kz\rangle =\langle\sigma w,Kz\rangle=\langle z,K\sigma
w\rangle.$
\end{proof}

\begin{corollary}\label{c1}
If $x,y\in V_{\mathbf{R}}$, then $\langle Hx,y\rangle $ is real
and $\langle Ax,y\rangle=\langle Ay,x\rangle$.
\end{corollary}

\begin{proof}
Clearly $\langle \cdot ,\cdot \rangle $ is real on $V_{\mathbb{R}}\times V_{%
\mathbb{R}}$. Since $HV_{\mathbb{R}}\subseteq V_{\mathbb{R}}$, we see $%
\langle Hx,y\rangle $ is real. Next, $ \langle Ax,y\rangle=
\langle Hx,y\rangle+ \langle Kx,y\rangle$. The term $\langle
Hx,y\rangle$ equals $\langle Hy,x\rangle$ because $\langle
Hx,y\rangle$ is real and $H$ is self-adjoint. On the other hand,
$\langle Kx,y\rangle=\langle Ky,x\rangle$ by Lemma \ref{l5}. So
$\langle Ax,y\rangle=\langle Ay,x\rangle$.
\end{proof}

\begin{lemma}
Define $m:V_{\mathbb{R}}\times V\rightarrow \mathbb{C}$ by $%
m(x,z)=e^{\langle Hz,x\rangle}e^{\langle K\bar {z},x\rangle}
e^{-\langle Ax,x\rangle/2}$. Then $m$ is a multiplier. Moreover,
if $T_{x}F(z):=m(x,z)F(z-x)$, then $x\mapsto T_x$ is a
representation of the abelian group $V_{\mathbb{R}}$ on
$\mathbf{F}_{A}$. It is unitary if $KV_{\mathbb{R}}\subseteq
V_{\mathbb{R}}$.
\end{lemma}

\begin{proof}
We first show $m$ is a multiplier:
\begin{align*}
m(x,&z)m(y,z-x)=e^{\langle Hz,x\rangle}e^{\langle K\bar
{z},x\rangle} e^{-\langle Ax,x\rangle /2}e^{ \langle
H(z-x),y\rangle}e^{ \langle
K(\bar {z}-x),y\rangle}e^{-\langle Ay,y\rangle /2}\\
&=e^{\langle Hz,x+y\rangle}e^{\langle K\bar
{z},x+y\rangle}e^{-\langle Hx,y\rangle}e^{-\langle
Kx,y\rangle}e^{-\langle Ax,x\rangle /2}e^{
-\langle Ay,y\rangle /2}\\
&=e^{\langle Hz,x+y\rangle}e^{\langle K\bar
{z},x+y\rangle}e^{-\langle Ax,y\rangle /2-\langle Ax,x\rangle /2-\langle Ay,y\rangle /2}\\
&=e^{\langle Hz,x+y\rangle}e^{\langle K\bar
{z},x+y\rangle}e^{-\langle
A(x+y),x+y\rangle /2}\\
&=m(x+y,z).
\end{align*}

Since $m$ is a multiplier, we have $T_xT_y=T_{x+y}$. For each $T_x$ to
be unitary, we need $|m(x,z)|=e^{(Az,x)-(Ax,x)/2}$. But
\begin{equation*}
|m(x,z)|=e^{(Hz,x)}e^{(K\bar {z},x)}e^{-(Ax,x)/2}=
e^{(Az,x)- (Ax,x)/2}e^{(K\bar{z}-Kz,x)}.
\end{equation*}
Thus $T_x$ is unitary for all $x$ if and only if the real part of
every vector $K\bar {z}-Kz$ is 0. Since $\bar {z}- z$ runs over
$iV_{\mathbb{R}}$ as $z$ runs over $V$, $T_x$ is unitary for all
$x$ if and only if $K(iV_{\mathbb{R}})\subset iV_{\mathbb{R}}$,
which is equivalent to $K(V_{\mathbb{R}})\subset V_{\mathbb{R}}$.

\end{proof}

Notice that $\det_{V}H=\left( \det H\right) ^{2}$. To simplify some
calculations later on we define $c:=(2\pi )^{-n/4}\left( \frac{\det_{V}A}{%
\det H}\right) ^{1/4}$. We remark for further reference:

\begin{lemma}
\label{lcons} $c_{A}^{-2}c^{2}=\frac{\sqrt{\det H}}{(2\pi )^{n/2}}$ and $%
c^{-1}\frac{\sqrt{\det (H)}}{\pi ^{n/2}}=\left( \frac{2}{\pi }\right) ^{n/4}\frac{%
\left( \det H\right) ^{3/4}}{\left( \det_{V}A\right) ^{1/4}}$.
\end{lemma}

Let $D(x)=c\,m(x,0)=c\,e^{-\langle Ax,x\rangle/2}$ and define $R:\mathbf{F}%
_{A}\rightarrow C^{\infty }(V_{\mathbb{R}})$ by $RF(x):=D(x)F(x)$. Extending
the bilinear form $x\mapsto \langle Ax,x\rangle$ to a complex bilinear form
 $\langle z,z\rangle_{A}$
on $V$ shows that $D$ has a holomorphic extension to $V$.

\begin{lemma}
The restriction map $R$ intertwines the action of $V_{\mathbb{R}}$ on $%
\mathbf{F}_{A}$ and the left regular action $L$ on functions on
$V_{\mathbb{R}}$.
\end{lemma}

\begin{proof}
We have
\begin{eqnarray*}
R(T_{y}F)(x) &=&c\,m(x,0)T_{y}F(x) \\
&=&c\,m(x,0)m(y,x)F(x-y) \\
&=&c\,m(x,0)m(-y,-x)F(x-y) \\
&=&c\,m(x-y,0)F(x-y) \\
&=&L_{y}RF(x).
\end{eqnarray*}
\end{proof}

\section{The Generalized Segal--Bargmann Transform}

As for the classical space,  $R$ is a densely defined, closed
operator.  It also has dense image in $L^{2}(V_{\mathbb{R}})$. To see this,
let $\left\{ h_{\alpha }\right\} _{\alpha }$ be the orthonormal
basis of $L^{2}(V_{\mathbb{R}})$ given by the Hermite functions. Then $%
\left\{ \det (A)^{\frac14}h_{\alpha }(\sqrt{A}x)\right\} _{\alpha }$ is an orthonormal
basis of $L^{2}(V_{\mathbb{R}})$ which is contained in the image of $R$. It
follows again that $R$ has an adjoint and
\begin{equation*}
R^{\ast }h(z)=(R^{\ast }h,K_{A,z})=(h,RK_{A,z})
\end{equation*}
where $K_{A,z}(w)=K_{A}(w,z)=c_{A}^{-2}e^{\frac{1}{2}\overline{\langle
Kw,w\rangle }}e^{\langle Hw,z\rangle }e^{\frac{1}{2}\langle Kz,z\rangle }$.
Thus
\begin{align*}
R^{\ast }h(z)& =c\,\int h(x)e^{-\langle Ax,x\rangle /2}\overline{K_{A}(x,z)}%
\,dx \\
& =c_{A}^{-2}c\int h(x)e^{-\langle Ax,x\rangle /2}e^{\frac{1}{2}\overline{%
\langle Kz,z\rangle }}e^{\langle z,Hx\rangle }e^{\frac{1}{2}\langle
Kx,x\rangle }\,dx \\
& =c_{A}^{-2}c\,e^{\frac{1}{2}\overline{\langle Kz,z\rangle }}\int
h(x)e^{-\langle Hx,x\rangle /2}e^{-\langle Kx,x\rangle /2}e^{\langle
z,Hx\rangle }e^{\frac{1}{2}\langle Kx,x\rangle }\,dx \\
& =c_{A}^{-1}c\,e^{\frac{1}{2}\overline{\langle Kz,z\rangle }}\int
h(x)e^{-\langle x,Hx\rangle /2}e^{\langle z,Hx\rangle }\,dx \\
& =c_{A}^{-2}c\,e^{\frac{1}{2}\overline{\langle Kz,z\rangle }}e^{\frac{1}{2}%
\langle z,H\bar{z}\rangle }\int h(x)e^{-(\langle z,H\bar{z}\rangle -\langle
z,Hx\rangle -\langle x,H\bar{z}\rangle +\langle x,Hx\rangle )/2}\,dx \\
& =c_{A}^{-2}c\,e^{\frac{1}{2}\overline{\langle Kz,z\rangle }}e^{\frac{1}{2}%
\langle z,H\bar{z}\rangle }\int h(x)e^{-\langle z-x,H(\bar{z}-\bar{x}%
)\rangle /2}\,dx
\end{align*}
for $\langle z,Hx\rangle =\langle \overline{Hx},\bar{z}\rangle =\langle Hx,%
\bar{z}\rangle =\langle x,H\bar{z}\rangle $ and $\langle z,Hx\rangle
=\langle z,H\bar{x}\rangle $. Thus we finally arrive at
\begin{equation}
R^{\ast }h(z)=c_{A}^{-2}c\,\,e^{\frac{1}{2}\langle z,H\bar{z}+Kz\rangle }e^{-%
\frac{1}{2}\langle x,H\bar{x}\rangle }\ast h(z).  \label{eq5.1}
\end{equation}
Let $P:V_{\mathbb{R}}\rightarrow V_{\mathbb{R}}$ be positive. Define
$\phi _{P}(x)=\sqrt{\det (P)}(2\pi )^{-n/2}e^{-||\sqrt{P}x||^{2}/2}$.
For $t>0$, let $P(t)=P/t$.

\begin{lemma}
\label{l9}Let the notation be as above.\ Then $0<t\mapsto \phi _{P(t)}$ is a
convolution semigroup, i.e., $\phi _{P(t+s)}=\phi _{P(t)}\ast \phi _{P(s)}$.
\end{lemma}

\begin{proof}
This follows by change of parameters $y=\sqrt{P}x$ from the fact that $\phi
_{Id(t)}(x)=(2\pi t )^{-n/2}e^{-||x||^{2}/2t}$ is a convolution semigroup.
\end{proof}

We define a unitary operator $W$ on $L^2(V_{\mathbb R})$ by
\[Wf(x)=e^{i\mathrm{Im}\langle x,Kx\rangle}f(x)=e^{i\mathrm{
Im}\langle x,Ax\rangle}f(x).\]
We know $W=I$ if $KV_{\mathbb R}\subseteq V_{\mathbb R}$ and this occurs
if $A$ leaves $V_{\mathbb R}$ invariant.

\begin{lemma}\label{l10}
Let $h$ be in the domain of definition of $R^{\ast }$. Then
$RR^{\ast }h=W(\phi _{H}*h)$.
\end{lemma}

\begin{proof}
We notice first that $c_{A}^{-2}c^{2}=(2\pi )^{-n/2}\sqrt{\det H\,}$ by
Lemma \ref{lcons}. From (\ref{eq5.1}) we then get
\begin{align*}
RR^{\ast}h(x)&=c\,e^{-\frac 12\langle Ax,x\rangle}R^{\ast}h(x)\notag\\
&=c_A^{-2}c^2\,e^{-\frac 12\langle Ax,x\rangle}e^{\frac 12%
\langle x,H\bar {x}+Kx\rangle}e^{-\frac 12\langle y,H\bar {y}\rangle}
\ast h(x)\notag\\
&=(2\pi )^{-n/2}\sqrt {\det(H)}\,e^{-\frac 12\langle Ax,x\rangle}
e^{%
\frac 12\langle x,Ax\rangle}e^{-\frac 12\langle y,H\bar {y}\rangle}
\ast h(x)\notag\\
&=(2\pi )^{-n/2}\sqrt {\det(H)}\,e^{i\text{\rm Im}\langle x,Ax\rangle}
\int e^{-\frac 12(y,Hy)}h(x-y)\,dy.\notag\\
&=(2\pi )^{-n/2}\sqrt {\det(H)}e^{i\text{\rm Im}\langle x,Ax\rangle}
\int e^{-\frac {||\sqrt Hy||^2}2%
}h(x-y)\,dy\notag\\
&=W(\phi_H\ast h)(x)\notag
\end{align*}
\end{proof}

Lemma \ref{l9} and Lemma \ref{l10} leads to the following corollary:
\begin{corollary}
Suppose $AV_{\mathbb R}\subseteq V_{\mathbb R}$. Then
\begin{equation*}
|R^{\ast}|h(x)=\phi_{H(1/2)}\ast h(x)=\frac {\sqrt {\det(H)}}{\pi^{n/2}}
\int_{V_{\mathbb{R}}}e^{-||\sqrt {H}y||^2}h(x-y)\,dy.
\end{equation*}
\end{corollary}

\begin{theorem}[The Segal--Bargmann Transform]
Suppose $A$ leaves $V_{\mathbb R}$ invariant.
Then the operator $U_A:L^2(V_{\mathbb{R}})\rightarrow\mathbb{F}_A$ defined
by
\begin{equation*}
U_{A}f(z)=\left( \frac{2}{\pi }\right) ^{n/4}\frac{\left( \det H\right)
^{3/4}}{\left( \det_{V}A\right) ^{1/4}}e^{\frac{1}{2}\left( \langle Hz,\bar{z%
}\rangle +\langle z,Kz\rangle \right) }\int e^{\langle H(z-y),\bar{z}%
-y\rangle }f(y)\,dy\,.
\end{equation*}
is a unitary isomorphism. The map $U_A$ is called the \textbf{%
generalized Segal--Bargmann transform.}
\end{theorem}

\begin{proof}
By polarization we can write $R^{\ast }=U\,|R^{\ast }|$ where $U:L^{2}(V_{%
\mathbb{R}})\rightarrow \mathbf{F}_{A}$ is a  unitary isomorphism.
Taking adjoints gives $|R^{\ast }|\,U^{\ast }=R$. Hence
$RU=|R^{\ast }|$. Hence
\begin{align*}
c\,m(x)Uh(x)& =RUh(x) \\
& =(|R^{\ast }|\,h)(x) \\
& =\frac{\sqrt{\det (H)}}{\pi ^{n/2}}\int_{V_{\mathbf{R}}}e^{-||\sqrt{H}%
y||^{2}}h(x-y)\,dy.
\end{align*}
Since $m(x)=e^{-\frac{1}{2}(\langle x,Hx\rangle +\langle x,Kx\rangle )}$, we
have using
Lemma \ref{lcons}:
\begin{equation*}
Uf(x)=\left( \frac{2}{\pi }\right) ^{n/4}\frac{\left( \det H\right) ^{3/4}}{%
(\det_{V}A)^{1/4}}\,e^{\frac{1}{2}\left( \langle x,Hx\rangle +\langle
x,Kx\rangle \right) }\int e^{(x-y,H(x-y))}f(y)\,dy.
\end{equation*}
By holomorphicity, this implies
\begin{equation*}
Uf(z)=\left( \frac{2}{\pi }\right) ^{n/4}\,\frac{\left( \det H\right) ^{3/4}%
}{(\det_{V}A)^{1/4}}e^{\frac{1}{2}\left( \langle Hz,\bar{z}\rangle +\langle
z,Kz\rangle \right) }\int e^{\langle H(z-y),\bar{z}-y\rangle }f(y)\,dy
\end{equation*}
is the Bargmann--Segal transform.
\end{proof}

\section{The Gaussian Formulation}

In infinite dimensions, there is no useful notion of Lebesgue
measure but Gaussian measure does make sense. So, with a view to
extension to infinite dimensions, we will recast our generalized
Segal-Bargmann transform using Gaussian measure instead of
Lebesgue measure as the background measure on $V_{\mathbf{R}}$. Of
course, we have already defined the Fock space $\mathbf{F}_A$
using Gaussian measure.

As before, $V$ is a finite-dimensional complex vector space with
Hermitian inner-product $\langle\cdot,\cdot\rangle$, and $A:V\to
V$ is a {\it real}--linear map which is {\it symmetric,
positive-definite} with respect to the real inner-product
$(\cdot,\cdot)=\text{%
\textrm{Re}}\langle \cdot,\cdot\rangle$, i.e. $(Az,z)>0$ for all
$z\in V$ except $z=0$. We assume, furthermore, that there is a
real subspace $ V_{\mathbb{R}}$ for which
$V=V_{\mathbb{R}}+iV_{\mathbb{R}}$, the inner-product
$\langle\cdot,\cdot\rangle$ is real-valued on $V_{\mathbb{R}}$ and
$A(V_{\mathbb{R}})\subset V_{\mathbb{R}}$. As usual, $A$ is the
sum
$$A=H+K$$
where $H=(A-iAi)/2$ is complex-linear on $V$ and $K=(A+iAi)/2$ is
complex-conjugate-linear. The real subspaces $V_{\mathbb{R}}$ and
$iV_{\mathbb{R}}$ are $(\cdot,\cdot)$-orthogonal because for any
$x,y\in V_{\mathbb{R}}$ we have $(x,iy)=\textrm{Re}\langle x,
iy\rangle=-\textrm{Re}(i\langle x,  y\rangle)$, since $\langle x,
y\rangle$ is real, by hypothesis. Since $A$ preserves
$V_{\mathbb{R}}$ and is symmetric, it also preserves the
orthogonal complement $iV_{\mathbb{R}}$. Thus $A$ has the block
diagonal form
\begin{equation}
A=\left[\begin{matrix} R&0\\0&T\end{matrix}\right]=d(X,Y)
\end{equation}
Here, and henceforth, we   use the notation $d(X,Y)$ to mean the
real-linear map $ V\to V$ given by $a\mapsto Xa$ and $ia\mapsto
iYa$ for all $a\in V_{\mathbb{R}}$, where $X, Y$ are real-linear
operators on $V_{\mathbb{R}}$.   Note that $d(X,Y)$ is
complex-linear if and only if $X=Y$ and is
complex-conjugate-linear if and only if $Y=-X$. The operator
$d(X,X)$ is the unique complex-linear map $V\to V$ which restricts
to $X$ on $V_{\mathbb{R}}$, and we will denote it by $X_V$:
\begin{equation}\label{eq:defXV}
X_V=\left[\begin{matrix} X&0\\0&X\end{matrix}\right]
\end{equation}

The hypothesis that $A$ is symmetric and positive-definite (by
which we mean $A>0$, not just $A\geq 0$) means that $R$ and $T$
are symmetric, positive definite on $V_{\mathbb{R}}$.
Consequently, the real-linear operator $S$ on $V_{\mathbb{R}}$
given by
\begin{equation}\label{eq:defS}
S=2(R^{-1}+T^{-1})^{-1}
\end{equation}
is also symmetric, positive-definite.

The operators $H$ and $K$ on $V$ are given by
\begin{equation}\label{eq:HKRT}H= \frac{1}{2}(R_V+T_V) ,\quad K =
d\left(\frac{1}{2}(R-T), \frac{1}{2}(T-R)\right) \end{equation}
Using the conjugation map $$\sigma:V\to V:a+ib\mapsto
a-ib\qquad\mbox{for $a,b\in V_{\mathbb{R}}$ }$$ we can also write
$K$ as
\begin{equation}\label{eq:KsigmaRT}
K=\frac{1}{2}(R_V-T_V)\sigma \end{equation}
 Now consider the
holomorphic functions $\rho_T$ and $\rho_S$ on $V$ given by
\begin{eqnarray}\label{eq:defrhoP}
\rho_T(z)=\frac{(\det
T)^{1/2}}{(2\pi)^{n/2}}e^{-\frac{1}{2}\langle T_Vz,
\bar{z}\rangle}\\
\rho_S(z)=\frac{(\det
S)^{1/2}}{(2\pi)^{n/2}}e^{-\frac{1}{2}\langle S_Vz,
\bar{z}\rangle}
\end{eqnarray}
where $n=\dim V_{\mathbb{R}}$. Restricted to $V_{\mathbb{R}}$,
these are density functions for Gaussian probability measures.

The \textbf{Segal-Bargmann transform} in this setting is given by
the map
\begin{equation}\label{eq:defSA}
S_A:L^2(V_{\mathbb{R}},\rho_S(x)dx)\to \mathbf{F}_{A} : f\mapsto
S_Af
\end{equation}
where
\begin{equation}\label{eq:defSAfz}
S_Af(z)=\int_{V_{\mathbb{R}}}f(x)\rho_T(z-x)\,dx=\int_{V_{\mathbb{R}}}f(x)c(x,z)\rho_S(x)\,dx
\end{equation}
where the generalized ``coherent state'' function $c$ is
specified,  for $x\in V_{\mathbb{R}}$ and $z\in V$,   by
\begin{equation}\label{eq:defLxw}
c(x,z)=\frac{\rho_T(x-z)}{\rho_S(x)}
\end{equation}
It is possible to take (\ref{eq:defSAfz}) as the starting point,
with $f\in L^2(V_{\mathbb{R}},\rho_S(x)dx)$ and prove that:
(i)$S_Af(z)$ is well-defined, (ii) $S_Af$ is in $\mathbf{F}_{A}$,
(iii)$S_A$ is a unitary isomorphism onto $\mathbf{F}_{A}$.
However, we shall not work out everything in this approach since
we have essentially proven all this in the preceding sections.
Full details of a direct approach would be obtained by
generalizing the procedure used in \cite{AS99}. In the present
discussion we shall work out only some of the properties of $S_A$.

\begin{lemma}\label{lem:Lxw}: Let $w,z\in V$. Then:
\begin{itemize}
\item[(i)] The function $x\mapsto c(x,z)$ belongs to
$L^2(V_{\mathbb{R}},\rho_S(x)dx)$, thereby ensuring that the
integral (\ref{eq:defSAfz}) defining $S_Af(z)$ is well-defined;
\item[(ii)] The $S_A$--transform of $c(\cdot,w)$ is $K_A(\cdot,\bar{w})$:
\begin{equation}\label{eq:SAl}
[S_Ac(\cdot,w)](z)=K_A(z,\bar{w})
\end{equation}
and so, in particular,
\begin{equation}\label{eq:KAintegralexpression}
K_A(z,w)=\int_{V_{\mathbb{R}}}\frac{\rho_T(x-z)\rho_T(x-{\bar w}
)}{\rho_S(x)}\,dx
\end{equation}
\item[(iii)] The transform $S_A$ preserves inner--products on the
linear span of the functions $c(\cdot,w)$:
 $$\langle
c(\cdot,w),c(\cdot,z)\rangle_{L^2(V_{\mathbb{R}},\rho_S(x)dx)}=
K_A(w,z)=\langle K_A(\cdot,{\bar w} ),K_A(\cdot,{\bar z}
)\rangle_{\mathbf{F}_{A}}$$
\end{itemize}
\end{lemma}
\begin{proof} (i) is equivalent to finiteness of
$\int_{V_{\mathbb{R}}}\frac{|\rho_T(x-z)|^2}{\rho_S(x)}\,dx$,
which is equivalent to positivity of the operator   $2T-S$. To see
that $2T-S$ is positive observe that
\begin{eqnarray}
2T-S & = &\nonumber
2T[(R^{-1}+T^{-1})-T^{-1}](R^{-1}+T^{-1})^{-1}\nonumber\\
&=&   2TR^{-1}(R^{-1}+T^{-1})^{-1}=TR^{-1}S\label{eq:2tminussshort}\nonumber\\
&=& 2(T^{-1}+T^{-1}RT^{-1})^{-1}\label{eq:2tminuss}
\end{eqnarray}
and in this last line $T^{-1}>0$ (being the inverse of $T>0$) and
$(T^{-1}RT^{-1}x,x)=(RT^{-1}x,T^{-1}x)\geq 0$ by positivity of
$R$. Thus $2T-S$ is positive, being twice the inverse of the
positive operator $T^{-1}+T^{-1}RT^{-1}$.

(ii) is the result of a lengthy calculation which, despite an
unpromising start, leads from complicated expressions to simple
ones. To avoid writing a lot of complex conjugates we shall use
the symmetric complex bilinear pairing $v\cdot w=vw=\langle
v,\bar{w}\rangle$ for $v,w\in V$, writing $v^2$ for $vv$. More
seriously, we shall denote the complex-linear operator  $T_V$
which restricts to $T$ on $V_{\mathbb{R}}$ simply by $T$. It is
readily checked that $T$ continues to be symmetric in the sense
that $Tv\cdot w=v\cdot Tw$ for all $v,w\in V$. We start with
\begin{eqnarray*}
a &\stackrel{\rm def}{=}& [S_Ac(\cdot,w)](z)\\
&=&\int_{V_{\mathbb{R}}}\frac{\rho_T(x-w)}{\rho_S(x)}\rho_T(z-x)\,dx\\
&=& (2\pi)^{-n/2}\frac{\det T}{(\det S)^{1/2}
}\int_{V_{\mathbb{R}}} e^{-\frac{1}{2} \left[T(x-w)\cdot
(x-w)+ T(x-z)\cdot (x-z)-Sx\cdot x\right]}\,dx\\
&=&(2\pi)^{-n/2}\frac{\det T}{(\det S)^{1/2}
}\int_{V_{\mathbb{R}}} e^{-\frac{1}{2}\left[(2T-S)x\cdot
x-2Tx\cdot (w+z)+Tw\cdot w+Tz\cdot z\right]}\,dx
\end{eqnarray*}
Recall from the proof of (i) that $2T-S>0$. For notational
simplicity let $L=(2T-S)^{1/2}$ and $M=L^{-1}T$. Then
\begin{eqnarray*}
a&=&(2\pi)^{-n/2}\frac{\det T}{(\det S)^{1/2}
}\int_{V_{\mathbb{R}}}e^{-\frac{1}{2}\left(Lx-M(w+z)\right)^2}\,dx\,e^{-\frac{1}{2}
\left[Tw\cdot w+Tz\cdot z-M(w+z)\cdot M(w+z)\right]}\\
&=& \frac{\det T}{(\det S)^{1/2}(\det L) }e^{-\frac{1}{2}
\left[Tw\cdot w+Tz\cdot z-M(w+z)\cdot M(w+z)\right]}
\end{eqnarray*}
To simplify the last exponent observe that
\begin{eqnarray*}Tw\cdot w-Mw\cdot Mw &=& Tw\cdot w -Tw\cdot
L^{-2}Tw\\
&=&Tw\cdot w -Tw\cdot (2T-S)^{-1}Tw\\
&=& Tw\cdot w-\frac{1}{2}Tw\cdot (T^{-1}+T^{-1}RT^{-1})Tw\qquad\mbox{using (\ref{eq:2tminuss})}\\
&=& Tw\cdot w-\frac{1}{2}Tw\cdot (w+T^{-1}Rw)\\
&=&\frac{1}{2}(Tw\cdot w-Rw\cdot w)\\
&=& - \langle K{\bar w},{\bar w}\rangle \qquad\mbox{by
(\ref{eq:KsigmaRT})}
\end{eqnarray*}
The same holds with $z$ in place of $w$. For the ``cross term'' we
have
\begin{eqnarray*} Mw\cdot Mz &=& Tw\cdot L^{-2}Tz\\
&=& Tw\cdot (2T-S)^{-1}Tz\\
&=& \frac{1}{2}Tw\cdot (T^{-1}+T^{-1}RT^{-1})Tz\\
&=& \frac{1}{2}(Tw\cdot z+w\cdot Rz)\\
&=& 2  w\cdot Hz
\end{eqnarray*}
Putting everything together we have
$$ [S_Ac(\cdot,w)](z)=\frac{\det T}{(\det S)^{1/2}(\det L)
}e^{\frac{1}{2}\langle K{\bar w},{\bar w}\rangle} e^{ \langle
Hw,\bar{z}\rangle
 }e^{\frac{1}{2}\langle K{\bar z},{\bar z}\rangle}$$
 In Lemma \ref{lemma:cAdet} below we prove that
 $$\frac{\det T}{(\det S)^{1/2}(\det L)
}=\left( \frac{\det_{V}(A)}{\det_{V}(H)}\right) ^{-1/2}=c_A^{-2}$$
 So
$$ [S_Ac(\cdot,w)](z)=K_A(w,\bar{z})$$

For (iii), we have first:
$$
\langle c(\cdot,w),c(\cdot,z)\rangle_{L^2(\rho_S(x)dx)} =
[S_Ac(\cdot,w)](\bar{z})=K_A({\bar z},\bar{w})=K_A(w,z)$$ The
second equality in (iii) follows from the fact that $K_A$ is a
reproducing kernel. \end{proof}

\section{The evaluation map and   determinant relations  }

Recall the reproducing kernel
$$K_A(z,w)=c_A^{-2}e^{\frac{1}{2}\langle z,Kz\rangle+\frac{1}{2}\langle Kw,w\rangle
+\langle Hz,w\rangle}$$ where
$$c_A^{-2}=\left(\frac{\det_VH}{\det_VA}\right)^2$$
Being a reproducing kernel for $\mathbb{F}_A$ means
\begin{equation}\label{eq:fKAw}
f(w)=\bigl(f,K_A(\cdot,w)\bigr)=\pi^{-n}(\det
A)^{1/2}\int_Vf(z)K_A(w,z)\ |dz|
\end{equation}
where $|dz|=dxdy$ signifies integration with respect to Lebesgue
measure on the real inner-product space $V$. Thus we have
\begin{proposition} For any $z\in V$, evaluation map
$$ \delta_z:\mathbb{F}_A\to\mathbb{C}:f\mapsto f(z)$$ is
bounded linear functional with norm
\begin{equation}\label{eq:normdeltaz}
\|\delta_z\|=K_A(z,z)^{1/2}=c_A^{-1}e^{(Az,z)}
\end{equation}

\end{proposition}
\begin{proof} We have
\begin{equation}\label{eq:deltazf}
|\delta_zf|=|f(z)|=|(f,K_A(\cdot,z))|\leq
\|f\|_{\mathbb{F}_A}K_A(z,z)^{1/2} \end{equation} because, again by
the reproducing kernel property we have
$$\|K_A(\cdot,z)\|_{\mathbb{F}_A}^2=\bigl(K_A(\cdot,z),
K_A(\cdot,z)\bigr)_{\mathbb{F}_A}=K_A(z,z)$$ This last calculation
also shows that the inequality in (\ref{eq:deltazf}) is an
equality of $f=K_A(\cdot,z)$ and thereby shows that $\|\delta_z\|$
is actually equal to $K_A(z,z)^{1/2}$. The latter is readily
checked to be equal to $c_A^{-1}e^{(Az,z)}$.
\end{proof}
Next we make two observations about the constant $c_A$, the first
of which has already been used.

\begin{lemma}\label{lemma:cAdet}For the constant $c_A$ we have
$$c_A^{-2}=\left(\frac{\det_VH}{\det_VA}\right)^2=\frac{\det
T}{(\det S)^{1/2}\det L}$$ where, as before, $L=(2T-S)^{1/2}$ and
$S=2(R^{-1}+T^{-1})^{-1}$.
\end{lemma}
\begin{proof}
Recall from (\ref{eq:2tminussshort}) that $2T-S=TR^{-1}S$. Note
also that
$$S^{-1}=\frac{1}{2}(R^{-1}+T^{-1})=R^{-1}\frac{R+T}{2}T{-1}=R^{-1}(H|V_{\mathbb{R}})T^{-1}$$
So
\begin{eqnarray*}
\left(\frac{\det_VA}{\det_VH}\right)^{1/2}\frac{\det T}{(\det
S)^{1/2}\det L} &=&\frac{(\det R)^{1/2}(\det T)^{1/2}}{\det
S^{-1}\det R\det T}\frac{\det T}{(\det S)^{1/2}\det T^{1/2}\det
R^{-1/2}\det S^{1/2}}\\
&=&1 \end{eqnarray*} which implies the dersired result.
\end{proof}

Next we prove a determinant relation which implies $c_A\geq 1$:
\begin{lemma} If $R$ and $T$ are positive definite $n\times n$ matrices
(symmetric if real) then
\begin{equation}\label{eq:detinequality}
\sqrt{\det R \det T}\leq \det \left(\frac{R+T}{2} \right)
\end{equation}
with equality if and only if $R=T$.

\end{lemma}
\begin{proof}
Note first that the matrix \begin{equation} D\stackrel{\rm def}{=}
R^{-1/2}TR^{-1/2} \end{equation}
 is positive definite because
$(R^{-1/2}TR^{-1/2}x,x)=(TR^{-1/2}x,R^{-1/2}x)\geq 0$ since
$T>0$, with equality if and only if $R^{_1/2}x=0$ if and only if
$x=0$. So $D=(R^{-1/2}TR^{-1/2})^{1/4}$ makes sense and is also
positive definite (and is symmetric if we are working with reals).
We have then
\begin{eqnarray*}
\frac{\det R \det T}{\left(\det \frac{R+T}{2}\right)^2}&=&
\frac{\det R\det (R^{1/2}D^4R^{1/2})}{\left[\det
R^{1/2}\left(\frac{1+D^4}{2}\right)R^{1/2}\right]^2} \\
&=& \left[\det\left(\frac{D^2+D^{-2}}{2}\right)\right]^{-2}\\
&=&
\left[\det\left\{I+\left(\frac{1}{\sqrt{2}}D-\frac{1}{\sqrt{2}}
D^{-1}\right)^2\right\}\right]^{-2}
\end{eqnarray*}
To summarize:
\begin{equation}\label{eq:determinantidentity}
\frac{\det R \det T}{\left(\det \frac{R+T}{2}\right)^2}=
\left[\det\left\{I+\left(\frac{1}{\sqrt{2}}D-\frac{1}{\sqrt{2}}
D^{-1}\right)^2\right\}\right]^{-2}
\end{equation}
where $ D=(R^{-1/2}TR^{-1/2})^{1/4}$. Diagonalizing $D$ makes it
apparent that this last term is $\leq 1$ with equality if and only
if $D=D^{-1}$, which is equivalent to $D^4=I$ which holds if and
only if $R=T$.
\end{proof}

As consequence we have for $c_A$:
$$c_A=\left(\frac{\det_VA}{\det_VH}\right)^{1/4}=\left(
\frac{\det R \det
T}{\left(\det\frac{R+T}{2}\right)^2}\right)^{1/4} = \left(\frac{
\sqrt{\det R \det T}}{ \det\frac{R+T}{2}  }\right)^{1/2}  $$ and
so
\begin{equation}\label{eq:cAminusttwo}
c_A^{-2}=\frac{ \det\frac{R+T}{2} }{ \sqrt{\det R \det T}}\geq 1
\end{equation}
 with equality if and only if $R=T$.

When extending this theory to infinite-dimensions we have to note
that in order to retain a meaningful notion of evaluation
$\delta_z: f\mapsto f(z)$, the constant $c_A^{-1}$ which appears
in the norm $\|\delta_z\|$ given in (\ref{eq:normdeltaz}) must be
finite. The expression for $c_A^{-2}$ obtained from
(\ref{eq:determinantidentity}) gives a more explicit condition on
$R$ and $T$ for this finiteness to hold.

If $R$ and $T$ are both scalar operators, say $R=rI$ and $T=tT$,
then (\ref{eq:cAminusttwo}) shows that $c_A^{-1}$ equals
$[(r+t)/(2\sqrt{rt})]^{n/2}$ which is bounded as
$n\nearrow\infty$ if and only if $r=t$. This observation was made
in \cite{AS99}.

\section{Remarks on extension to infinite dimensions}

The Gaussian formulation permits extension to the
infinite-dimensional situation, at least with some conditions
placed on $A$ so as to make such an extension reasonable. Suppose
then that $V$ is an infinite-dimensional separable complex Hilbert
space, $V_{\mathbb{R}}$ a real subspace on which the inner-product
is real-valued, and $A:V\to V$ a bounded symmetric,
positive-definite real-linear operator carrying $V_{\mathbb{R}}$
into itself. The operators $R$, $T$, $S$, $H$ and $K$ are defined
as before. Assume  that $R$ and $T$ commute and that there is an
orthonormal basis $e_1, e_2,...$ of $V_{\mathbb{R}}$ consisting of
simultaneous eigenvectors  of $R$ and $T$ (greater generality may
be possible but we discuss only this case) . Let $V_n$ be the
complex linear span of $e_1,...,e_n$, and $V_{n,\mathbb{R}}$ the
real linear span of $e_1,...,e_n$. Then $A$ restricts to an
operator $A_n$ on $V_n$, and we have similarly restrictions $H_n,
K_n$ on $V_n$ and $R_n, T_n, S_n$ on $V_{n,\mathbb{R}}$. The
unitary transform $S_A$ may be obtained as a limit of the
finite-dimensional transforms $S_{A_n}$.

The Gaussian kernels $\rho_S$ and $\rho_T$ do not make sense
anymore, and nor does the coherent state $c$, but the Gaussian
measures $d\gamma_S(x)= \rho_S(x)dx$ and $\mu_A$ do have
meaningful analogs. There is a probability space
$V'_{\mathbb{R}}$, with a $\sigma$--algebra ${\mathcal F} $ on
which there is a measure $\gamma_S$, and there is a linear map
$V_{\mathbb{R}}\to L^2(V'_{\mathbb{R}},\gamma_A):x\mapsto
G(x)=(x,\cdot)$, such that the $\sigma$--algebra ${\mathcal F} $
is generated by the random variables $G(x)$, and each $G(x)$ is
(real) Gaussian with mean $0$ and variance $(S^{-1}x,x)$.
Similarly, there is probability space $V'$, with a
$\sigma$--algebra ${\mathcal F}_1 $ on which there is a measure
$\mu_A$, and there is a real-linear map $V\to
L^2(V',\mu_A):z\mapsto G_1(z)=(z,\cdot)$, such that the
$\sigma$--algebra ${\mathcal F}_1 $ is generated by the random
variables $G_1(z)$, and each $G_1(z)$ is (real) Gaussian with mean
$0$ and variance $\frac{1}{2}(A^{-1}z,z)$. Then for each $z\in V$,
written as $z=a+ib$ with $a,b\in V_{\mathbb{R}}$, we have the
complex-valued random variable on $V'$ given by
$${\tilde z}=G_1(a)+iG_1(b)$$
Suppose $g$ is a holomorphic function of $n$ complex variables
such that
$$\int_V|g({\tilde e}_1,...,{\tilde
e}_n)|^2\,d\mu_A<\infty.$$
 Define ${\mathbb{F}}_A$ to be the
closed linear span of all functions of the type $g({\tilde
e}_1,...,{\tilde e}_n)$ in $L^2(\mu_A)$ for all $n\geq 1$. We may
then define $S_A$ of a function $f\bigl(G(e_1),...,G(e_n)\bigr)$
to be $(S_{A_n} f)({\tilde e}_1,...,{\tilde e}_n)$, and then
extend $S_A$ be continuity to all of $L^2(\gamma_S)$. In writing
$(S_{A_n} f)$ we have identified $V_n$ with ${\mathbb{C}}^n$ and $
V_{n,\mathbb{R}}$ with $\mathbb{R}^n$ using the  basis
$e_1,...,e_n$.

A potentially significant application of the  infinite-dimensional
case would be to situations where $V_{\mathbb{R}}$ is a {\it path
space} and $A$ is arises from a suitable {\it differential
operator}. For the ``classical case'' where $R=T=tI$ for some
$t>0$, this leads to the Hall transform \cite{H94} for Lie groups
as well as the path-space version on Lie groups considered in
\cite{HS98}.


\end{document}